\documentclass[preprint]{elsarticle}

\usepackage{hyperref}
\usepackage{lineno}
\usepackage{amsmath}
\usepackage{amsfonts}
\usepackage{amssymb}
\usepackage[margin=1in]{geometry}
\usepackage{color}
\usepackage[all]{xy}
\usepackage{graphicx}

\journal{Discrete Mathematics}

\newtheorem{theorem}{Theorem}

\newtheorem{axiom}[theorem]{Axiom}

\newtheorem{definition}[theorem]{Definition}
\newdefinition{rmk}{Remark}
\newdefinition{notation}{Notation}
\newdefinition{example}[theorem]{Example}
\newproof{proof}{Proof}

\newproof{pot}{Proof of Theorem \ref{thm2}}

\newcommand{\GG}{\mathfrak{G}}
\newcommand{\II}{\mathfrak{I}}

\begin{document}

\begin{frontmatter}

\title{About a counterexample on contractible transformations of graphs}

%
\author{Mart\'in-Eduardo Fr\'ias-Armenta}
\ead{martineduardofrias@gmail.com}
\ead{eduardo.frias@unison.mx}
%
%
\address{Departamento de Matm\'aticas, Universidad de Sonora, M\'exico.}

\begin{abstract}
In this paper we show that most of the results in \cite{Ivashchenko:contractible} are false or the proofs are wrong.
\end{abstract}


\end{frontmatter}


\section{Introduction}

In graph theory, several reductions that leave the homology invariant have been studied. In \cite{Ivashchenko:homology,Ivashchenko:contractible}, A. Ivashchenko shows a family of graphs constructed from $K(1)$ by contractible transformations (as in Definition \ref{TI}), and he proves that such transformations do not change the homology groups of graphs. He started the study of these transformations because these are used in the theory of molecular spaces and digital topology. Modern references are \cite{Boutry,Han,Espinoza}. Particulary 
 \cite{Espinoza} contributes with 
\begin{enumerate}
\item an aplication to topological data analisys,
\item they have proved that the Ivashchenko's  contractible graphs are collapsibles 
\item and we can also see as a collorary from the main result of \cite{Ivashchenko:homology} i.e. that the homology do not change by the Ivashchenko's transformations.
\end{enumerate}
We show in this paper that most of the results in \cite{Ivashchenko:contractible} are false or the proofs are wrong.

\begin{figure}[h] 
\begin{center}
 \includegraphics[scale=0.3]{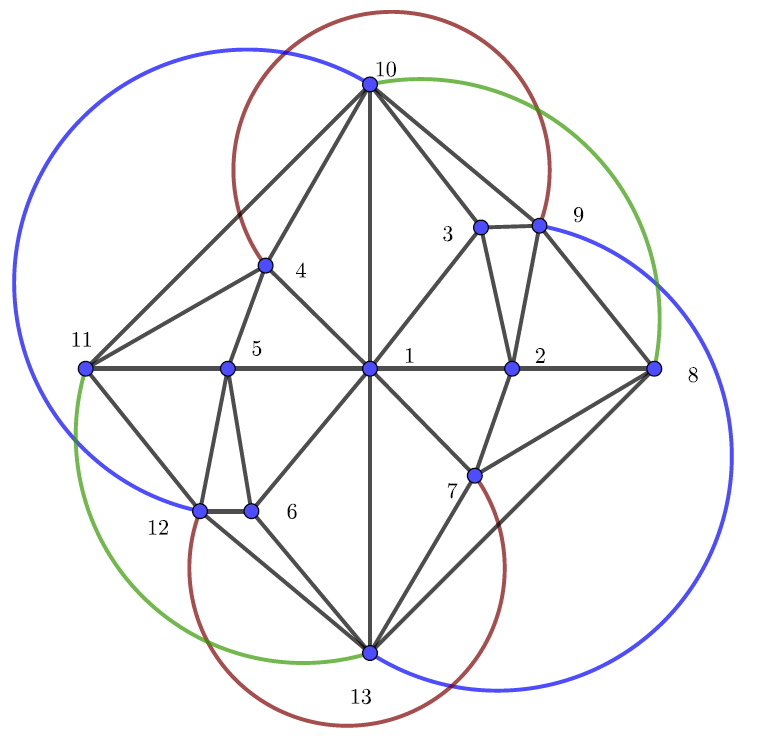}
\caption{Heart graph. We can see in the heart graph that the axiom  3.4 of \cite{Ivashchenko:contractible} is false. We can not put any edge in vertex $1$ as the axiom 3.4 in that paper states. 
}\label{figura}
\end{center}
\end{figure}


\section{Contractible transformations}
Let $G=(V(G), E(G)) \in \GG$ be a graph, and let $v \in V(G)$ be a vertex. We denote  $N_G(v) = \{ u \in V(G): \{u,v\} \in E(G)\}$ and $N_G(v,w) = N_G(v) \cap N_G(w)$. In addition, by abuse of notation we identify the graph with its set of vertices.

In \cite{Ivashchenko:homology} the next family of graphs was defined, and its elements are called contractible graphs.

\begin{definition}\label{TI} Let $\II \subset \GG$ be the family of graphs defined by
\begin{enumerate}
\item The trivial graph $K(1)$ is in $\II$.
\item Any graph of $\II$ can be obtained  from $K(1)$  by the following transformations.
\begin{enumerate}
\item[(I1)] Deleting a vertex $v$. A vertex $v$ of a graph $G$ can be deleted if $N_G(v)\in \II$.
\item[(I2)] Gluing  a vertex $v$. If a subgraph $G_1$ of the graph $G$ is in $\II$, then the vertex $v$ can be glued to the graph $G$ in such way that $N_G(v) = G_1$.
\item[(I3)] Deleting  an edge $\{v_1,v_2\}$. The edge $\{v_1,v_2\}$ of a graph $G$ can be deleted if $N_G(v_1,v_2) \in \II$.
\item[(I4)] \label{gluedge} Gluing an edge $\{v_1,v_2\}$. Let two vertices $v_1$ and $v_2$ of a graph $G$ be nonadjacent. The edge $\{v_1,v_2\}$ can be glued if $N_G(v_1,v_2) \in \II$.
\end{enumerate} 
\end{enumerate}
If $G$ belongs to $\II$, then $G$ is called a contractible graph.
\end{definition}

The transformations (\textit{I1})-(\textit{I4}) were referred in \cite{Ivashchenko:homology} as contractible transformations. The contractible transformations are used in molecular spaces, see \cite{Ivashchenko:homology} for more explanation. In addition, in \cite{Ivashchenko:homology} it was proved that contractible transformations do not change the homology groups of a graph, for any commutative group of coefficients $A$, so the elements of $\II$ have trivial groups of $A$-homology.

In \cite[Th. 4.9]{Ivashchenko:homology}, it was proved that the contractible transformation 
 does not change the homology groups of the graph $G$ when $N_G(v) \in \II$.
\section{Counterexample}

We cite textualy the Axiom 3.4 at \cite{Ivashchenko:contractible}

\begin{axiom}
  Suppose  that  $G$ is  a  contractible  graph,  and  a  vertex  $v$, $v\in    G$,  is  not  adjacent 
to  some  vertices  of  $G$.  Then  there  exists  a  nonadjacent  vertex  $u$, $u\in G$,  such  that  the 
subgraph  $O(vu)$  is  contractible. 

\end{axiom}
Where $O(vu)$ is the induce graph by $N(v)\cap N(u)$.\\
Ivashchenko claimed that the previous axiom is verified on small graphs and he did not intent to prove the generic case. 
But in the heart graph of figure \ref{figura}, we can see that vertex $1$ is not adjecent to $8$, $9$, $11$ y $12$, 
we can see that common neighboorhood of $1$ with each of those vertices is not contractible. So the Ivashchenko's axiom is false.
All the results of that paper are based in the axiom 3.4, and so the Theorems 3.5, 3.8, 3.9, 3.10 and Corollary 3  are clearly false,
 the heart graph in figure \ref{figura} shows this. For example theorem 3.5 establishes that any contractible graph has two contractible vertices and in figure \ref{figura} we clearly see that the heart graph is contractible and it does not have any contractible vertex.
 We think that theorem 3.7, 3.11 and 3.12 are true, but the proofs are incorrect because they use the axiom 3.4 or some of its consecuenses. 
The only theorem that is correctly proved it is 3.13.\\

The heart graph in figure \ref{figura} is the smallest graph visually pleasing that we found. We are wondering if there is  another with less vertices.

\section*{Acknowledgment}
I  thank Anton Dochthermann, Jes\'us Espinoza, Etiene Fieux and H\'ector Hern\'andez for their friendly and solidary support in the writing of this note.

\bibliographystyle{elsart-num}

\end{document}